\documentclass[12pt]{article}
\usepackage{amsmath,amssymb}
\usepackage[dvips]{graphicx,color}
\usepackage{a4wide}
\usepackage{natbib}
\usepackage{mathrsfs}
\usepackage{bm}
\begin{document}

\newcommand{\bqn}{\begin{eqnarray*}}
\newcommand{\eqn}{\end{eqnarray*}}
\newcommand{\bqa}{\begin{eqnarray}}
\newcommand{\eqa}{\end{eqnarray}}
\newcommand{\um}{{\underline{m}}}
\newcommand{\cC}{{\mathcal C}}
\newcommand{\bA}{\mathbf{A}}
\newcommand{\bB}{\mathbf{B}}
\newcommand{\bF}{\mathbf{F}}
\newcommand{\bG}{\mathbf{G}}

\newcommand{\bT}{\mathbf{T}}
\newcommand{\bX}{\mathbf{X}}
\newcommand{\bY}{\mathbf{Y}}

\newcommand{\bSi}{\boldsymbol{\Sigma}}
\newcommand{\conds}{Conditions [A]-[B]}
\newcommand{\cov}{\mathop{\mathrm{cov}}}
\newcommand{\CC}{{\mathbb{C}}}
\newcommand{\bgma}{{\mbox{\boldmath$\gamma$}}}
\newcommand{\gD}{\Delta}
\newcommand{\gd}{\delta}
\newcommand{\bbe}{{\bf e}}
\newcommand{\bbA}{{\bf A}}
\newcommand{\bbB}{{\bf B}}
\newcommand{\bbS}{{\bf S}}
\newcommand{\bbT}{{\bf T}}
\newcommand{\bbI}{{\bf I}}
\newcommand{\rtr}{{\rm tr}}
\newcommand{\rE}{{\rm E}}
\newcommand{\ointctrclockwise}{\oint}
\newcommand{\proof}{\noindent {\bf Proof.}}
\newcommand{\cU}{{\cal U}}
\newcommand{\rdd}{\textcolor{red}}
\newcommand{\eprf}{\hspace*{\fill}$\blacksquare$}
\numberwithin{equation}{section}  % number the equations by section

\newtheorem{thm}{Theorem}[section]
\newtheorem{lem}{Lemma}[section]
\newtheorem{rem}{Remark}[section]
\newtheorem{cor}{Corollary}[section]
\newtheorem{exm}{Example}[section]

\baselineskip 24pt

%%%%%%%%%%%%%%%%%%%%%%%%%%%%%%%%%%%%%%%%%%%%%%%%%%%%%%%%%%%%%%%%%%%%%%%%%%%%%%%%%%%%%%%%%%%%%%%%%%%%%%%%%%

\begin{center}
 {\bf\Large A Note on Central Limit Theorems for Linear Spectral
Statistics of Large Dimensional \emph{F}-matrix}\\
\end{center}
\vspace{.4cm} \centerline{Shurong Zheng and Zhidong Bai} \vspace{.4cm} \centerline{\it KLAS and School of Mathematics $\&$ Statistics, Northeast Normal
University, P.R.China} \vspace{.55cm}
 \centerline{\it E-mails: zhengsr@nenu.edu.cn; baizd@nenu.edu.cn}
\fontsize{9}{11.5pt plus.8pt
minus .6pt}\selectfont

%%%%%%
\begin{quotation}
\noindent {\it Abstract.} Sample covariance matrix and multivariate $F$-matrix play important roles in multivariate statistical analysis.
The central limit theorems {\sl (CLT)} of linear spectral statistics associated with these matrices were established in Bai and Silverstein (2004) and Zheng (2012) which received considerable attentions and have been applied to solve many large dimensional statistical problems. However, the sample covariance matrices used in these papers are not centralized  and there exist some questions about CLT's defined by the centralized sample covariance matrices. In this note, we shall provide some short complements on the CLT's in Bai and Silverstein (2004) and Zheng (2012),
 and show that the results in these two papers remain valid for the centralized sample covariance matrices, provided that the ratios of dimension $p$ to sample sizes $(n,n_1,n_2)$ are redefined as $p/(n-1)$ and $p/(n_i-1)$, $i=1,2$, respectively.

%\par\vspace{9pt}\noindent
%{\it AMS 2000 subject classifications.} Primary 15A52, 60F05; secondary 62H10.
\par
\end{quotation}\par

\section{Introduction}
%Sample covariance matrix and F matrix play an important role in statistics (see Anderson (2003); Bai (1999)).
%Many statistical inference will involve sample covariance matrix and F matrix.
%For example, seeking confidence intervals
%of multiple correlation coefficients and specifying graphical structures
%will use large dimensional sample covariance matrices.
%testing the equality of covariance matrices and canonical correlation analysis will involve
%large dimensional F matrix (see Bai and Saranadasa(1996); Bai (2009)).

Let $\{X_{jk}, j,k=1,2,\cdots\}$ and $\{Y_{jk}, j,k=1,2,\cdots\}$ be
two independent double arrays of independent random variables, either both real or both complex. In the sequel, we use $A^*$ to denote a complex conjugate transpose of a vector or matrix $\bA$. For $p>
1$, $n>1$ and $N>1$, we define $\bX =({\bX}_{1},\cdots,{\bX}_{ n})$ and $\bY =({\bY}_{1},\cdots,{\bY}_{ N})$ with  column vectors ${\bf X}_{j} = (X_{j1},...,X_{jp})', 1\leq j \leq n$, and ${\bf Y}_{k} = (Y_{k1},...,Y_{kp})'$, $1\leq k \leq N$. Let $\bT_p$ be a $p\times p$ non-negative definite (nnd) matrix. There exists a unique nnd matrix $\bT_p^{1/2}$ such  that $\bT_p=(\bT_p^{1/2})^2$. Then, $({\bf T}_p^{1/2}{\bf X}_{1},\cdots,{\bf T}_p^{1/2}{\bf X}_{
n})$ and $({\bf T}_p^{1/2}{\bf Y}_{1},\cdots,{\bf T}_p^{1/2}{\bf Y}_{ N})$ can be considered as two independent samples of sizes $n$ and $N$, respectively, drawn from
a $p$-dimensional population with population covariance matrix ${\bf T}_p$.

It is well known that the sample covariance matrices for $\bT_p^{1/2}\bX$ and $\bT_p^{1/2}\bY$ are often defined as
\bqa
{\bf S}_x&=&\frac1{n-1}\left(\sum_{i=1}^{n}{\bf T}_p^{1/2}{\bf X}_i{\bf X}_i^*{\bf T}_p^{1/2}-n{\bf T}_p^{1/2}\bar{\bf X}\bar{\bf X}^*{\bf T}_p^{1/2}\right),\label{eq11}\\
{\bf S}_y&=&\frac1{N-1}\left(\sum_{i=1}^{N}{\bf T}_p^{1/2}{\bf Y}_i{\bf Y}_i^*{\bf T}_p^{1/2}-N{\bf T}_p^{1/2}\bar{\bf Y}\bar{\bf Y}^*{\bf T}_p^{1/2}\right)\label{eq12},
\eqa
respectively, where $\bar {\bf X}=\frac1{n}
\sum\limits_{i=1}^n{\bf X}_i$ and $\bar {\bf Y}=\frac1{N}\sum\limits_{i=1}^N{\bf Y}_i$. The multivariate $F$-matrix\footnote{To guarantee that the definition makes sense, we need to assume that $p<N$ and ${\bf T}_p$ is positive definite. Because the eigenvalues of ${\bf F}$ are independent of ${\bf T}_p$,
we may assume ${\bf T}_p$ is an identity matrix.} is then defined as
\bqa
{\bf F}={\bf S}_x{\bf S}_y^{-1}. \label{eq13}
\eqa
Notice that the matrices defined in (\ref{eq11})--(\ref{eq13}) are transformation invariant, we will call them {\em centralized} sample covariance matrices and multivariate $F$-matrix, respectively.

Due to Corollary A.41 and Theorem A.43 of Bai and Silverstein (2009), in the literature of random matrix theory, the sample covariance matrices are usually {\em simplified} as
\bqa
{\bf B}_x=\frac1{n}\sum_{i=1}^{n}{\bf T}_p^{1/2}{\bf X}_i{\bf X}_i^*{\bf T}_p^{1/2} \ \mbox{and}\
{\bf B}_y=\frac1{N}\sum_{i=1}^{N}{\bf T}_p^{1/2}{\bf Y}_i{\bf Y}_i^*{\bf T}_p^{1/2}
\eqa
and the multivariate $F$-matrix is {\em simplified} as ${\bf G}={\bf B}_x{\bf B}_y^{-1}$.

Bai and Silverstein (2004) considered
the central limit theorem (CLT) of the linear spectral statistics (LSS) of the simplified sample covariance matrix
${\bf B}_x$ and provided the explicit expressions of asymptotic means and covariance functions for $\bB_x$. Later, Zheng (2012) extended the work of Bai and Silverstein (2004)
to the case of the multivariate $F$-matrix ${\bf G}$ and obtained explicit expressions of the asymptotic means, variances, and covariances for $\bG$.

Examining the inequalities derived from  Corollary A.41 and Theorem A.43 of Bai and Silverstein (2009), one finds that the difference between the empirical spectral distributions ({\sl ESD}) of ${\bf S}_x$ and ${\bf B}_x$ is of the order $O(n^{-1})$. Hence, we conclude that
${\bf S}_x$ and ${\bf B}_x$ have the same limiting spectral distributions ({\sl LSD}). However, the scale normalizers in CLT's of LSS of random matrices ${\bf S}_x$ and ${\bf B}_x$ have the same order as $p$. Thus, it is expected that the asymptotic biases in the CLT's of LSS of ${\bf S}_x$ and ${\bf B}_x$
should have a little difference. Upon such a consideration,
Pan (2012) reconsidered the CLT of LSS of centralized sample covariance matrix ${\bf S}_x$.
To reduce the asymptotic bias, he added an additional term to that of Bai and Silverstein (2004), that is,
 \begin{equation}
 \frac{y}{2\pi i}\int g(z)\frac{\underline{m}_{y}(z)\int\frac{tdH(t)}{(1+t\underline{m}_{y}(z))^2}}{z\left(1-
 y\int\frac{\underline{m}_{y}^2(z)t^2dH(t)}{(1+t\underline{m}_{y}(z))^2}\right)}dz,
 \label{eq15}
 \end{equation}
 where $\underline{m}_{y}(z)=-\frac{1-y}{z}+y m_{y}(z)$, $m_{y}(z)$ is the Stieltjes transform of
 the LSD of ${\bf S}_x$, $H(t)$ is the LSD of ${\bf T}_p$
and $y_{n}=p/n\rightarrow y> 0$.

%That is,
%in computing the asymptotic mean, we have to compute an additional term which is not easy to compute in practice.
%Zheng (2012) obtained CLT of LSS of $F$-matrix
%$\begin{array}{lll}
%\bF:={\bf S}_1{\bf S}_2^{-1}
%\end{array}$
% where
% $$
% {\bf S}_1=\frac{1}{n} \sum\limits_{k=1}^{n}{\bf X}_{\cdot k} {\bf X}_{\cdot k}^{*},\quad
%{\bf S}_2=\frac{1}{N} \sum\limits_{k=1}^{N}{\bf Y}_{\cdot k} {\bf Y}_{\cdot k}^{*}.
% $$
% However, due to the unknown population mean vectors of ${\bf X}_{\cdot i}$ and ${\bf Y}_{\cdot i}$, people are willing to
% use unbiased sample covariance matrices
%$$\tilde{{\bf S}}_1=\frac{1}{n-1}\sum\limits_{i=1}^{n}({\bf X}_{\cdot i}-\bar{\bf X})
%({\bf X}_{\cdot i}-\bar{\bf X})^{*}\quad\mbox{and}\quad\tilde{{\bf S}}_2=\frac{1}{N-1}
%\sum\limits_{i=1}^{N}({\bf Y}_{\cdot i}-\bar{\bf Y})({\bf Y}_{\cdot i}-\bar{\bf Y})^{*}.$$
%Then a general F matrix is defined as follows
%$
%\tilde{\bF}:=\tilde{{\bf S}}_1\tilde{{\bf S}}_2^{-1}.
%$ This paper will obtain that the LSS of $\tilde{{\bf S}}={\bf T}_p^{1/2}\tilde{{\bf S}}_1{\bf T}_p^{1/2}$
%and LSS of ${\bf S}={\bf T}_p^{1/2}{\bf S}_1{\bf T}_p^{1/2}$ have the same CLT.
%Moreover, the LSS of $\tilde{{\bf F}}=\tilde{{\bf S}}_1\tilde{{\bf S}}_2^{-1}$
%and LSS ${\bf F}={\bf S}_1{\bf S}_{2}^{-1}$ have the same CLT. Then although
%population mean vector are unknown, unbiased sample
%covariance matrix and the corresponding F matrix still have the same LSD and CLTs as
%Bai and Silverstein (2004) and Zheng (2012).

It is well known that when the population is multivariate-normally distributed,  the centralized sample covariance matrix
 ${\bf S}_x$ has the same distribution as simplified covariance matrix ${\bf B}_x$ with sample size $n-1$ and population mean zero. This fact motivates that
 this phenomenon should be asymptotically true in the general case. In this note, we shall give short proofs to indicate that if the simplified sample
 covariance matrix ${\bf B}_x$ is replaced by centralized
sample covariance matrix ${\bf S}_x$, Bai and Silverstein (2004)'s result remains valid provided that
the ratio of dimension to sample size $y_n$ is replaced by $p/(n-1)$
(this is equivalent to $c_n=n/(N-1)$ in  Bai and Silverstein (2004)). This result is equivalent to but much simpler than
that of Pan (2012) in both expressions and proof. Moreover, we shall prove that if the simplified multivariate $F$-matrix ${\bf G}$
is replaced by the centralized ${\bf F}$,
the results of Zheng (2012) remain valid provided the ratios of
dimensions to sample sizes,  $y_{n1}$ and $y_{n2}$, are replaced by $p/(n-1)$ and $p/(N-1)$.

  The remainder of this note is arranged as follows: Section 2 states the main theorems and the proof of Theorem \ref{Newthm2}. Section 3
  gives the proof of Theorem \ref{Newthm1}. The technical lemmas and their proofs will be postponed to Section \ref{sec4}.

\section{Main Results}
As mentioned in the previous section, the centralized covariance matrix will have the same LSD as that of the corresponding simplified covariance matrix. In this note, we shall prove the following theorems.
\begin{thm}\label{Newthm1}
Assume that
\par
(a) For each $p$, $\{X_{ij}, i\leq p,j\leq n\}$ are independent random variables with $EX_{ij}=0$, $E|X_{11}|^2=1$,
and satisfying
\begin{equation}
      \frac{1}{np}\sum_{j=1}^p \sum_{k=1}^{n}
      E|X_{jk}|^4
      \mathbf{1}_{\{|X_{jk}|\geq \eta \sqrt{n}\}}
      \rightarrow 0,
      \qquad \mbox{for any fixed } \eta>0.
      \label{eq21}
    \end{equation}
    Note that the random variables may be allowed to depend on $p$, but we suppress this dependence from the notation for brevity.
\par
(b) We assume $E|X_{ij}|^4=3$ for the real case, and $E|X_{ij}|^4=2$ and $EX_{ij}^2=0$ for the complex case.
\par
(c) $y_{n}=p/n\rightarrow y$, and
\par
(d) ${\bf T}_p$ is a $p\times p$ non-random nnd  Hermitian  matrix with bounded spectral norm in $p$, and its ESD
$H_p\stackrel{D}{\rightarrow}H$ where $H$ is a proper probability distribution.
\par
Let $f$  be an analytic function on an open region in the complex plane which covers the support of LSD of ${\bf S}_x$ with the origin excluded.

Then
\par
(i) the random variables
\begin{equation}\label{1.5}
X_p(f)=p\int f(x)d\left( F^{{\bf S}_x}-F^{\{y_{n-1},H_p\}}(x)\right),
\end{equation}
form a tight sequence in $p$, where $F^{{\bf S}_x}$
is the ESD of centralized sample covariance matrix ${\bf S}_x$, $F^{\{y,H\}}$ is the LSD of ${\bf S}_x$ whose
LSD's Stieltjes transform
$m_y(z)$ satisfies $\underline m_y(z)=ym_y(z)-(1-y)/z$ and $\underline m_y(z)$ is the unique solution to the equation
\bqa
%\underline{m}^{\{y_{n},H_{p}\}}(z)=-\frac{1-y_{n}}{z}+y_{n}m^{\{y_{n},H_{p}\}}(z)
%\quad\mbox{and}\quad
z=-\frac{1}{\underline{m}_y}+{y}\int\frac{t}{1+t\underline{m}_y(z)}
dH(t).
\label{eq23}
\eqa
in the upper half complex plane for each $z\in \mathbb{C}^+=\{z:\Im(z)>0\}$.
\par
(ii) The random variables in (\ref{1.5}) converges weakly to Gaussian variables $X_{f}$
with the same means
%$$
%EX_{f_j}=-\frac{\kappa-1}{2\pi i}\int f_j(z)\frac{{y}\int\underline{m}_{y}^3(z)t^2(1+t\underline{m}_{y}(z))^{-3}dH(t)}
%{\left(1-{y}\int\underline{m}_{y}^2(z)t^2(1+t\underline{m}_{y}(z))^{-2}dH(t)\right)^2}dz
%$$
and covariance functions as given in Theorem 1.1 of Bai and Silverstein (2004).
\end{thm}
The proof of Theorem 2.1 is postponed to Section \ref{sec3}.

As for the CLT of LSS of ${\bf F}$ matrix, we have the following theorem.

\begin{thm}\label{Newthm2}

  Assume that
  \begin{enumerate}
  \item the two arrays $\{X_{jk}, j\leq p, k\leq n\}$ and $\{Y_{jk}, j\leq p, k\leq N\}$ satisfy for any fixed $\eta>0$,
  \begin{equation}
      \frac{1}{np}\sum_{j=1}^p \sum_{k=1}^{n}
      E|X_{jk}|^4
      \mathbf{1}_{\{|X_{jk}|\geq \eta \sqrt{n}\}}
      \rightarrow 0,
      \qquad
      \frac{1}{Np}\sum_{j=1}^p \sum_{k=1}^{N}
      E|Y_{jk}|^4
      \mathbf{1}_{\{|Y_{jk}|\geq \eta \sqrt{N}\}}
      \rightarrow 0.
      \label{eq24}
    \end{equation}
  \item For all $j,k$,
    $|EX_{jk}^4|=\beta_x+1+\kappa$,
    $|EY_{jk}^4|=\beta_y+1+\kappa$.
    If both ${\bf X}$ and ${\bf Y}$ are complex valued, then
    $EX^2_{jk}=EY^2_{jk}=0$. Moreover,
     $y_n=p/n\to y_1>0$ and $y_N=p/N\to y_2\in(0,1)$.
  \end{enumerate}
  Let  $f$  be an analytic function  in an open region
  of the complex plane containing the interval
  $\left[\frac{(1-h)^2}{(1-y_2)^2},\frac{(1+h)^2}{(1-y_2)^2}\right]$, the support of the continuous part of the LSD
  $F_{\bf y}$ of {\bf F}-matrix, $h=\sqrt{y_1+y_2-y_1y_2}$ and ${\bf y}=(y_1,y_2)$.

  Then, as $p\to\infty$, the random variables
 $$
W_p(f)=p\int f(x)d\left(F^{\bf F}(x)-F_{(y_{n-1},y_{N-1})}(x)\right)
$$
  converges weakly to
  Gaussian variables $\{W_{f}\}$ which have the same means
and covariance functions as given in Zheng (2012)
,where $F^{{\bf F}}(x)$ is the ESD of centralized $F$-matrix ${\bf F}$
and $F_{(y_1,y_2)}(x)$ is the LSD defined by (2.4)
of Zheng (2012).
\label{Thmnew2}
\end{thm}
Proof. As mentioned in Section 1, we may assume ${\bf T}_p$ to be an identity matrix. Split our proofs into two steps by writing
$$
tr({\bf F}-z{\bf I}_p)^{-1}-pm_{(y_{n-1},y_{N-1})}(z)=
\left[tr({\bf S}_x{\bf S}_y^{-1}-z{\bf I}_p)^{-1}-pm^{(y_{n-1},F^{{\bf S}^{-1}_y})}(z)\right]
+p\left[m^{(y_{n-1},F^{{\bf S}^{-1}_y})}(z)-m_{(y_{n-1},y_{N-1})}(z)\right]
$$
where $F^{{\bf S}_y^{-1}}(t)$ and $F^{{\bf S}_y}(t)$ are the ESDs of ${\bf S}_y^{-1}$ and ${\bf S}_y$, $m_{(y_1,y_2)}$ is the Stieltjes transform of the LSD of F matrix,
$$
\underline{m}^{\{y_{n-1},F^{{\bf S}_y^{-1}}\}}=-\frac{1-y_{n-1}}{z}+
y_{n-1}m^{\{y_{n-1},F^{{\bf S}_y^{-1}}\}}(z)
$$
\bqa
z&=&-\frac{1}{\underline{m}^{\{y_{n-1},F^{{\bf S}_y^{-1}}\}}}
+y_{n-1}\int\frac{t}{1+t\underline{m}^{\{y_{n-1},F^{{\bf S}_y^{-1}}\}}}dF^{{\bf S}_y^{-1}}(t)\nonumber\\
&=&-\frac{1}{\underline{m}^{\{y_{n-1},F^{{\bf S}_y^{-1}}\}}}
+y_{n-1}\int\frac{1}{t+\underline{m}^{\{y_{n-1},F^{{\bf S}_y^{-1}}\}}}dF^{{\bf S}_y}(t).\label{eq25}
\eqa
Step 1. Given ${\bf S}_y$ , in the proof of Theorem \ref{Newthm1}, we have proved that the process
$tr({\bf S}_x{\bf S}_y^{-1}-z{\bf I}_p)^{-1}-pm^{\{y_{n-1},F^{{\bf S}_y^{-1}}\}}(z)$
weakly tends to a Gaussian process on the contour with mean and covariance function as given in (6.29) and (6.30) of Zheng (2012).

Step 2. By (\ref{eq25}) and the truth of
\bqa
z&=&-\frac{1}{\underline{m}_{\{y_{n-1},y_{N-1}\}}}
+y_{n-1}\int\frac{1}{t+\underline{m}_{\{y_{n-1},y_{N-1}\}}}dF_{y_{N-1}}(t),\label{eq26}
\eqa
where $F_{y_{N-1}}$ is the M-P law with ratio of dimension to sample size $y_{N-1}$. Subtracting both sides of (\ref{eq25}) from those of
(\ref{eq26}), we obtain
\bqa\label{F}
&&{p\cdot\left[m^{\{y_{n-1},F^{{\bf S}^{-1}_y}\}}(z)-m_{\{y_{n-1},y_{N-1}\}}(z)\right]}
={N\cdot\left[\underline{m}^{\{y_{n-1},F^{{\bf S}^{-1}_y}\}}(z)-
\underline{m}_{\{y_{n-1},y_{N-1}\}}(z)\right]}\nonumber \\
&=&{-y_{n-1}\underline{m}_{\{y_{n-},y_{N-1}\}}\underline{m}^{\{y_{{n-1}},
F^{{\bf S}^{-1}_y}\}}\frac{tr\left({\bf S}_{y}+\underline{m}_{\{y_{n-1},y_{N-1}\}}{\bf I}_p\right)^{-1}
-pm_{y_{N-1}}(-\underline{m}_{\{y_{n-1},y_{N-1}\}})}
{1-y_{n-1}\cdot\int\frac{\underline{m}_{\{y_{n-1},y_{N-1}\}}
\cdot\underline{m}^{\{y_{n-1},F^{{\bf S}^{-1}_y}\}}dF_{N-1}(t)}
{\left(t+\underline{m}_{\{y_{n-1},y_{N-1}\}}\right)\cdot
\left(t+\underline{m}^{\{y_{n-1},F^{{\bf S}^{-1}_x}\}}\right)}}}\nonumber\\
&=&-y_{n-1}\underline{m}_{\{y_{n-1},y_{N-1}\}}\underline{m}^{\{y_{{n-1}},
F^{{\bf S}^{-1}_y}\}}\frac{p[m_{N-1}(-\underline{m}_{\{y_{n-1},y_{N-1}\}})
-m_{y_{N-1}}(-\underline{m}_{\{y_{n-1},y_{N-1}\}})]}
{1-y_{n-1}\cdot\int\frac{\underline{m}_{\{y_{n-1},y_{N-1}\}}\cdot\underline{m}^{\{y_{n-1},
F^{{\bf S}^{-1}_y}\}}dF_{N-1}(t)}
{\left(t+\underline{m}_{\{y_{n-1},y_{N-1}\}}\right)\cdot \left(t+\underline{m}^{\{y_{n-1},
F^{{\bf S}^{-1}_y}\}}\right)}}\nonumber\\
&=&\frac{-y_{n-1}\underline{m}_{\{y_{n-1},y_{N-1}\}}\underline{m}^{\{y_{{n-1}},
F^{{\bf S}^{-1}_y}\}}\cdot
p[m_{N-1}(-\underline{m}_{\{y_{n-1},y_{N-1}\}})
-m_{y_{N-1}}(-\underline{m}_{\{y_{ n-1},y_{N-1}\}})]}{1-y_{n-1}\cdot\int\frac{\underline{m}_{\{y_{n-1},y_{N-1}\}}\cdot\underline{m}^{\{y_{n-1},
F^{{\bf S}^{-1}_y}\}}dF_{N-1}(t)}
{\left(t+\underline{m}_{\{y_{n-1},y_{N-1}\}}\right)\cdot \left(t+\underline{m}^{\{y_{n-1},
F^{{\bf S}^{-1}_y}\}}\right)}}
%&=&\frac{-y_{n-1}\underline{m}_{\{y_{n-1},y_{N-1}\}}\underline{m}^{\{y_{{n-1}},
%F^{{\bf S}^{-1}_y}\}}\cdot
%p[m_{N}(-\underline{m}_{\{y_{n-1},y_{N-1}\}})-m_{y_{N}}(-\underline{m}_{\{y_{n-1},y_{N-1}\}})]}
%{1-y_{n-1}\int\frac{\underline{m}_{\{y_{n-1},y_{N-1}\}}\cdot\underline{m}^{\{y_{n-1},
%F^{{\bf S}^{-1}_y}\}}dF_{N-1}(t)}
%{\left(t+\underline{m}_{\{y_{n-1},y_{N-1}\}}\right)\left(t+\underline{m}^{\{y_{n-1},
%F^{{\bf S}^{-1}_y}\}}\right)}}
%+o_{p}(1)
\eqa
which weakly tends to a Gaussian process on the contour with mean and covariance function as given in
(6.33) and (6.34) of Zheng (2012) where $m_{N-1}$ is the Stieltjes transform of ESD $F_{N-1}(x)$ of ${\bf S}_y$, $z=-\frac{1}{\underline{m}_{y_{N-1}}}+\frac{y_{N-1}}{1+\underline{m}_{y_{N-1}}}$ and
$\underline{m}_{y_{N-1}}(z)=-\frac{1-y_{N-1}}{z}+y m_{y_{N-1}}(z)$.
By Theorem 2.1 and (\ref{F}) we obtain that
$$
tr({\bf F}-z{\bf I}_p)^{-1}-pm_{\{ y_{n-1},y_{N-1}\}}(z)\quad\mbox{and}\quad
tr\left({\bf G}-z{\bf I}_p\right)^{-1}-pm_{\{y_{n},y_{N}\}}(z)
$$
have the same asymptotic distribution. Hence, the random variables
$$
\left(p\int f(x)d(F^{{\bf F}}(x)-F_{\{y_{n-1},y_{N-1}\}}(x))\right)
$$
converges weakly to Gaussian variables $W_{f}$ with the same means
and covariance functions as Zheng (2012).

Then the proof of Theorem 2.2 is completed. \eprf
\section{The Proof of Theorem \ref{Newthm1}}\label{sec3}

The condition (\ref{eq21}) allows us to truncate the random variables at $\eta_n\sqrt{n}$ and then renormalize them to have means zero and vairances 1, where $\eta_n\to 0$ slowly. Note that the 4th moments of the random variables may not be the same but they will be $\kappa+1+\beta_x+o(1)$ and for the complex case we have $EX_{ij}^2=o(n^{-1})$. The contour is defined similar to Bai and Silverstein (2004).
Define $\bgma_i=\frac{1}{\sqrt{n}}{\bf T}_p^{1/2}{\bf X}_{i}$. Then, we have
$$
{\bf S}_x=\frac{n}{n-1}\sum_{i=1}^n(\bgma_i-\bar \bgma)(\bgma_i-\bar \bgma)^*
=\sum_{i=1}^n\bgma_i\bgma_i^*-\frac{1}{n-1}\sum_{i\ne j}\bgma_i\bgma_j^*
={\bf B}_x-\gD
$$
where $\bar{\bgma}=\frac1n\sum\limits_{i=1}^n\bgma_i$ and $\Delta=\frac1{n-1}\sum_{j\ne k}\bgma_j\bgma_k^*$ .
%
%
%$$
%\begin{array}{lll}
%\bar{{\bf S}}&=&\frac{1}{n-1}\sum\limits_{i=1}^{n}({\bf T}_p^{1/2}{\bf X}_{\cdot i}-{\bf T}_p^{1/2}\bar{{\bf X}})
%({\bf T}_p^{1/2}{\bf X}_{\cdot i}-{\bf T}_p^{1/2}\bar{{\bf X}})^{*}\\
%&=&\frac{1}{n-1}\sum\limits_{i=1}^n{\bf T}_p^{1/2}{\bf X}_{\cdot i}{\bf X}_{\cdot i}^{*}{\bf T}_p^{1/2}-
%\frac{n}{n-1}{\bf T}_p^{1/2}\bar{{\bf X}}\bar{{\bf X}}^{*}{\bf T}_p^{1/2}\\
%&=&\frac{1}{n-1}\sum\limits_{i\not=k}{\bf T}_p^{1/2}{\bf X}_{\cdot i}{\bf X}_{\cdot i}^{*}{\bf T}_p^{1/2}
%+\frac{1}{n-1}{\bf T}_p^{1/2}{\bf X}_{\cdot k}{\bf X}_{\cdot k}^{*}{\bf T}_p^{1/2}
%-\frac{n}{n-1}{\bf T}_p^{1/2}\bar{{\bf X}}\bar{{\bf X}}^{*}{\bf T}_p^{1/2}
%\end{array}
%$$
%then we have
%$$
%\bar{{\bf S}}=\frac{1}{n}\sum\limits_{i}{\bf T}_p^{1/2}{\bf X}_{\cdot i}{\bf X}_{\cdot i}^{*}{\bf T}_p^{1/2}
%-\gD=\sum\limits_{i}\bgma_i\bgma_i^{*}-\gD={\bf S}-\gD
%$$
%where , ${\bf S}=\sum\limits_{i}\bgma_i\bgma_i^{*}$ and
%$\gD=\frac{1}{n(n-1)}\sum\limits_{k\not=j}{\bf T}_p^{1/2}{\bf X}_{\cdot k}{\bf X}_{\cdot j}^{*}{\bf T}_p^{1/2}
%=\frac{1}{n-1}\sum\limits_{k\not=j}\bm{\gamma}_k\bm{\gamma}_j^{*}.$
 Because
$$
\left({\bf S}_x-z{\bf I}\right)^{-1}=\left(\bbA-\gD\right)^{-1}=\bbA^{-1}
+\left(\bbA-\gD\right)^{-1}\gD\left(\bbA\right)^{-1}=\bbA^{-1}
+\bbA^{-1}\gD\bbA^{-1}+\bbA^{-1}(\gD\bbA^{-1})^2
+\left(\bbA-\gD\right)^{-1}\left(\gD\bbA^{-1}\right)^{3}
$$
where ${\bf A}(z)={\bf B}_x-z{\bf I}$, then we obtain
\begin{equation}\label{sn}
\begin{array}{lll}
&&p\left(\frac{1}{p}tr({{\bf B}_x}-z{\bf I})^{-1}-m_{n-1}^0(z)\right)\\
&=&p\left(\frac{1}{p}tr(\bbA-\gD)^{-1}-m_{n}^0(z)+m_{n}^0(z)-m_{n-1}^0(z)\right)\\
&=&p\left(\frac{1}{p}tr{\bf A}^{-1}(z)-m_{n}^0(z)\right)+p(m_{n}^0(z)-m_{n-1}^0(z))
+tr{\bf A}^{-2}(z)\gD\\
&&+tr{\bf A}^{-1}(z)(\gD{\bf A}^{-1}(z))^2+tr\left({\bf A}(z)-\gD\right)^{-1}(\gD{\bf A}^{-1}(z))^3\\
\end{array}
\end{equation}
where $\underline m_n^0(z)=-\frac{1-y_n}{z}+y_nm_n^0(z)$, $\underline m_n^0(z)$ and $\underline m_{n-1}^0(z)$ satisfy
\begin{equation}\label{eq1}
z=-\frac{1}{\underline{m}_{n}^{(0)}}+\frac{p}{n}\int\frac{t}{1+t\underline{m}_{n}^{(0)}}dH_p(t)
\end{equation}
\begin{equation}\label{eq2}
z=-\frac{1}{\underline{m}_{n-1}^{(0)}}+\frac{p}{n-1}\int\frac{t}{1+t\underline{m}_{n-1}^{(0)}}dH_p(t)
\end{equation}
\begin{equation}\label{eq12}
z=-\frac{1}{\underline{m}_{y}}+{y}\int\frac{t}{1+t\underline{m}_{y}}dH(t).
\end{equation}
By (\ref{sn}), Lemmas \ref{lem1} and \ref{thm2}, we have
\begin{equation}\label{end}
\rtr({{\bf B}_x}-z{\bf I})^{-1}-p\cdot m_{n-1}^0(z)=p\left(\frac{1}{p}
tr{\bf A}^{-1}(z)-m_{n}^0(z)\right)+o_p(1).
\end{equation}
So Theorem \ref{Newthm1} has the same asymptotic mean and covariance function as Theorem 1.1 of Bai and Silverstein
(2004).

\noindent
{\bf Tightness of $\rtr({{\bf S}_x}-z{\bf I})^{-1}-\rtr({{\bf B}_x}-z{\bf I})^{-1}$.} Using what has been proved in Bai and Silverstein (2004),
we only need to prove that there is an absolute constant $M$ such that for any $z_1,z_2\in {\cal C}$,
\bqa
&&\rE\Big|\rtr(\bbB_x-z_1\bbI)^{-1}(\bbB_x-z_2\bbI)^{-1}-\rtr (\bbS_x-z_1\bbI)^{-1}(\bbS_x-z_2\bbI)^{-1}\Big|^2\nonumber\\
&=&\rE\left|\sum_{i=1}^n\frac{(\lambda_i-\tilde\lambda_i)(\lambda_i+\tilde\lambda_i-2z_1z_2)}{(\lambda_i-z_1)
(\lambda-z_2)(\tilde\lambda_i-z_1)(\tilde\lambda-z_2)}\right|^2\nonumber\\
&\le& K\rE\left|\sum_{i=1}|\lambda_i-\tilde\lambda_i|\right|^2{\cal B}_n+o(1)\le M,\label{eq36}
\eqa
where $\{\lambda_i\}$ and $\{\tilde \lambda_i\}$ are the eigenvalues of $\bbB_x$ and $\bbS_x$, respectively, and arranged in descending order,
the event ${\cal B}_n$ is defined as $x_l+\epsilon <\tilde \lambda_p<\lambda_1<x_r-\epsilon$ such that
$$
{\rm P}({\cal B}_n)=o(n^{-3}).
$$
(for the justification of the definition ${\cal B}_n$, see Bai and Silverstein (1998)). The last step of (\ref{eq36}) follows from the fact that
$$
\sum_{i=1}|\lambda_i-\tilde\lambda_i|=\sum_{i=1}(\lambda_i-\tilde\lambda_i)\le \lambda_1-\tilde\lambda_p\le x_r
$$
by the interlacing theorem.

\noindent
{\bf The equi-continuity of $\rE\rtr({{\bf S}_x}-z{\bf I})^{-1}-pm^0_{n-1}(z)$} can be proved in a similar way to that for
the tightness of $\rtr({{\bf S}_x}-z{\bf I})^{-1}-\rE\rtr({{\bf B}_x}-z{\bf I})^{-1}$.

By now, the proof of Theorem 2.1 is completed. \eprf

\section{Technical Lemmas}\label{sec4}

\begin{lem}\label{lem1} Under assumptions of Theorem \ref{Newthm1}, for every $z\in{\mathbb C}^+$, we have
$$
p(m_{n}^{(0)}-m^{(0)}_{n-1})\rightarrow(1+z\underline{m}_{y})\cdot\frac{\underline{m}_{y}
+z\underline{m}'_{y}}{z\underline{m}_{y}}.
$$
\end{lem}
Proof. We have
\begin{equation}\label{eq3}
\underline{m}_{n}^{(0)}(z)=-\frac{n-p}{n}\cdot\frac{1}{z}+\frac{p}{n}m_{n}^0(z),\quad
\underline{m}_{n-1}^{(0)}(z)=-\frac{n-1-p}{n-1}\cdot\frac{1}{z}+\frac{p}{n}m_{n-1}^0(z)
\end{equation}
where $p/n\rightarrow {y}>0$. By (\ref{eq12}), we obtain
\begin{equation}\label{eq4}
\underline{m}'_{y}=\frac{1}{\frac{1}{\underline{m}_{y}^2}-{y}\int\frac{t^2}{(1+t\underline{m}_{y})^2}dH(t)},\quad
{y}\int\frac{t}{1+t\underline{m}_{y}}dH(t)=\frac{1+z\underline{m}_{y}}{\underline{m}_{y}}
\end{equation}
Using (\ref{eq1})-(\ref{eq2}), we obtain
$$
0=\frac{\underline{m}_{n}^{(0)}-\underline{m}_{n-1}^{(0)}}{\underline{m}_{n}^{(0)}\underline{m}_{n-1}^{(0)}}
-(\underline{m}_{n}^{(0)}-\underline{m}_{n-1}^{(0)})\frac{p}{n}\int\frac{t^2}{(1+t\underline{m}_{n}^{(0)})(
1+t\underline{m}_{n-1}^{(0)})}dH_p(t)-\frac{p}{n(n-1)}\int\frac{t}{1+t\underline{m}_{n-1}^{(0)}}dH_p(t),
$$
that is,
\begin{equation}\label{eq5}
n(\underline{m}_{n}^{(0)}-\underline{m}_{n-1}^{(0)})=\frac{\frac{p}{n-1}\int\frac{t}
{1+t\underline{m}_{n-1}^{(0)}}dH_p(t)}
{\frac{1}{\underline{m}_{n}^{(0)}\underline{m}_{n-1}^{(0)}}-
\frac{p}{n}\int\frac{t^2}{(1+t\underline{m}_{n}^{(0)})(1+t\underline{m}_{n-1}^{(0)})}dH_p(t)}\rightarrow
\frac{{y}\int\frac{t}{1+t\underline{m}_{y}}dH(t)}
{\frac{1}{\underline{m}_{y}^2}-
{y}\int\frac{t^2}{(1+t\underline{m}_{y})^2}dH(t)}.
\end{equation}
By (\ref{eq3}), (\ref{eq4}) and (\ref{eq5}), we have
\bqa
p(m_{n}^{(0)}-m^{(0)}_{n-1})&=&
n\underline{m}^{(0)}_{n}+\frac{n-p}{z}-\left((n-1)\underline{m}^{(0)}_{n-1}+\frac{n-1-p}z\right)\nonumber\\
&=&n(\underline{m}^{(0)}_{n}-\underline{m}^{(0)}_{n-1})+\underline m^{(0)}_{n-1}(z)+\frac{1}z\nonumber\\
&\rightarrow&
\frac{{y}\int\frac{t}{1+t\underline{m}_{y}}dH(t)}
{\frac{1}{\underline{m}_{y}^2}
-{y}\int\frac{t^2}{(1+t\underline{m}_{y})^2}dH(t)}+\frac{1+z\underline m(z)}{z}\nonumber\\
&=&\underline{m}_{y}^{*}\cdot\frac{1+z\underline{m}_{y}}{\underline{m}_{y}}+\frac{1+z\underline m(z)}{z}
=(1+z\underline{m}_{y})\cdot\frac{\underline{m}_{y}+z\underline{m}_{y}'}{z\underline{m}_{y}}.\label{eq6}
\eqa
Thus, we prove that Lemma \ref{lem1} holds.   \eprf

In the sequel, we shall use Vatali lemma frequently. Let
$$
\gD=\frac1n\sum_{j\ne k}\bgma_j\bgma_k^*.~(\mbox{It should be }\frac1{n-1}\sum_{j\ne k}\bgma_j\bgma_k^* \mbox{ but no harm to the limit.})
$$
We will derive the limit $\rtr(\bbA-\gD)^{-1}-\rtr(\bbA^{-1})$.

\begin{lem}\label{l1} After truncation and normalization, we have
$\rE|\bgma_k^*\bbA^{-1}\bgma_k-( 1+z\underline m_{y}(z))|^2\le K n^{-1}
$ for every $z\in{\mathbb C}^+$.
\end{lem}

\proof \ We have
$\bgma_k^*\bbA^{-1}\bgma_k=\bgma_k^*\bbA^{-1}_k\bgma_k\beta_k=1-\beta_k,$ where $\bbA_k=\bbA-\bgma_k\bgma_k^*$ and $\beta_k+(1+\bgma_k^*\bbA_k^{-1}\bgma_k)^{-1}$.
Therefore, By (1.15) and (2.17) of Bai and Silverstein (2004), we have
$\rE|\bgma_k^*\bbA^{-1}\bgma_k-g(z)|^2=\rE|\beta_k+z\underline{m}_{y}(z)|^2\le Kn^{-1}$
with $g(z)=1+z\underline{m}_{y}(z)$.

\begin{cor}After truncation and normalization, we have
$\rE\left|\bgma_k^*\bbA^{-2}\bgma_k-\frac{d}{dz}(1+z\underline{m}_{y}(z))\right|^2\le Kn^{-1}$
for every $z\in{\mathbb C}^+$.
\end{cor}
\proof\ By Cauchy integral formula, we have
$$
\bgma_k^*\bbA^{-2}\bgma_k=\frac1{2\pi i}\oint_{|\zeta-z|=v/2}\frac{\bgma_k^*\bbA^{-1}(\zeta)\bgma_k}{(\zeta-z)^2}d\zeta
\quad\mbox{and}\quad
g'(z)=\frac1{2\pi i}\oint_{|\zeta-z|=v/2}\frac{1+\zeta\underline{m}_{y}(\zeta)}{(\zeta-z)^2}d\zeta.
$$
Then $\rE\left|\bgma_k^*\bbA^{-2}\bgma_k-\frac{d}{dz}(1+z\underline{m}_{y}(z))\right|^2\le Kn^{-1}$ follows from Lemma \ref{l1}.

\begin{lem}\label{l2} For any subset $\cU$ of $\{1,2,\cdots,n\}$, after truncating and normalizing, we have
$\rE\left|\rtr{\bf A}^{-1}\gD\right|^2\le K n^{-1}$
for every $z\in{\mathbb C}^+$. Especially for every $z\in{\mathbb C}^+$,
$$
\rE|\rtr(\bbA^{-2}\gD)|^2=\rE\left|\frac1n\sum_{j\ne k}\bgma_j^*\bbA^{-2}\bgma_k\right|^2=O(n^{-1}).
$$
\end{lem}

\proof \ We have
$\rtr\bbA^{-1}\gD=\frac1n\sum_{j\ne k\in\cU}\bgma_j^*\bbA^{-1}\bgma_k=
\frac1n\sum_{j\ne k\in\cU}\bgma_j^*\bbA_{jk}^{-1}\bgma_k\beta_j\beta_{k(j)},$ where $\bbA_{jk}=\bbA_k-\bgma_j\bgma_j^*$ for $j\ne k$ and
$\beta_{k(j)}=(1+\bgma_k^*\bbA_{jk}^{-1}\bgma_k)^{-1}$. We will similarly define $\bbA_{ijk}$ and $\beta_{k(ij)}$ for later use.
Then we obtain
$$
\begin{array}{lll}
\rE|\rtr(\bbA^{-1}\gD)|^2&=&\rE\frac1n\sum_{j_1\ne k_1\in\cU}\bgma_{j_1}^*\bbA_{j_1k_1}^{-1}\bgma_{k_1}\beta_{j_1}\beta_{k_1(j_1)}
\frac1n\sum_{j_2\ne k_2\in\cU}\bgma_{k_2}^*\bbA_{j_2k_2}^{-1}\bgma_{j_2}\beta_{j_2}\beta_{k_2(j_2)}\\
&:=&\sum_{(2)}+\sum_{(3)}+\sum_{(4)},
\end{array}
$$
where the index $(\cdot)$ denotes the number of distinct integers in the set $\{j_1,k_1,j_2,k_2\}$. By the facts that $|\beta_{j}|\leq \frac{|z|}{\nu}$ and $\nu=\Im(z)$,
we have
$$
\begin{array}{lll}
\sum_{(2)}&\le &\frac{2|z|^4}{n^2v^4}\sum_{j\ne k\in\cU}\rE|\bgma_{j}^*\bbA_{jk}^{-1}\bgma_{k}|^2\\
&\leq&\frac{|z|^4}{\nu^4n^4}\sum_{j\ne k}\rE|\rtr(\bbT^*\bbA_{jk}^{-1}\bbT\bar \bbA_{jk}^{-1})\le
\frac{p}{n^2}\frac{|z|^4\|\bbT\|^2}{\nu^6}\le Kn^{-1}
\end{array}
$$
$$
\begin{array}{lll}
\sum_{(4)}&=&\frac{1}{n^2}\sum_{j_1\ne k_1\ne j_2\ne k_2\in\cU}
E\bgma_{j_1}^*\bbA_{j_1k_1}^{-1}\bgma_{k_1}\bgma_{j_2}^*\bbA_{j_2k_2}^{-1}\bgma_{k_2}
\beta_{j_1}\beta_{k_1(j_1)}\beta_{j_2}\beta_{k_2(j_2)}\\
\end{array}
$$
where
$$
\begin{array}{lll}
\bgma_{j_1}^*\bbA\bgma_{k_1}
&=&\beta_{j_1}\beta_{k_1(j_1)}\bgma_{j_1}^*\bbA_{j_1k_1}^{-1}\bgma_{k_1}
=\beta_{j_1}\beta_{k_1(j_1)}\Big[\bgma_{j_1}^*\bbA_{j_1k_1k_2}^{-1}\bgma_{k_1}-\beta_{k_2(j_1k_1)}\bgma_{j_1}^*\bbA_{j_1k_1k_2}^{-1}\bgma_{k_2}
\bgma_{k_2}^*\bbA_{j_1k_1k_2}^{-1}\bgma_{k_1}\Big]\\
&=&\beta_{j_1}\beta_{k_1(j_1)}\Big[\bgma_{j_1}^*\bbA_{j_1j_2k_1k_2}^{-1}\bgma_{k_1}-
\beta_{j_2(j_1k_1k_2)}\bgma_{j_1}^*\bbA_{j_1j_2k_1k_2}^{-1}\bgma_{j_2}
\bgma_{j_2}^*\bbA_{j_1j_2k_1k_2}^{-1}\bgma_{k_1}\\
&&-\beta_{k_2(j_1k_1)}\bgma_{j_1}^*\bbA_{j_1j_2k_1k_2}^{-1}\bgma_{k_2}
\bgma_{k_2}^*\bbA_{j_1j_2k_1k_2}^{-1}\bgma_{k_1}\\
&&+\beta_{k_2(j_1k_1)}\beta_{j_2(j_1k_1k_2)}
\bgma_{j_1}^*\bbA_{j_1j_2k_1k_2}^{-1}\bgma_{j_2}\bgma_{j_2}^*\bbA_{j_1j_2k_1k_2}^{-1}\bgma_{k_2}
\bgma_{k_2}^*\bbA_{j_1j_2k_1k_2}^{-1}\bgma_{k_1}\\
&&+\beta_{k_2(j_1k_1)}\beta_{j_2(j_1k_1k_2)}\bgma_{j_1}^*\bbA_{j_1j_2k_1k_2}^{-1}\bgma_{k_2}
\bgma_{k_2}^*\bbA_{j_1j_2k_1k_2}^{-1}\bgma_{j_2}\bgma_{j_2}^*\bbA_{j_1j_2k_1k_2}^{-1}\bgma_{k_1}\\
&&-\beta_{k_2(j_1k_1)}\beta^2_{j_2(j_1k_1k_2)}
\bgma_{j_1}^*\bbA_{j_1j_2k_1k_2}^{-1}\bgma_{j_2}\bgma_{j_2}^*\bbA_{j_1j_2k_1k_2}^{-1}\bgma_{k_2}
\bgma_{k_2}^*\bbA_{j_1j_2k_1k_2}^{-1}\bgma_{j_2}\bgma_{j_2}^*\bbA_{j_1j_2k_1k_2}^{-1}\bgma_{k_1}\Big]
\end{array}
$$
and
\bqa
\beta_j&=&{b}_j-\beta_j{b}_j\epsilon_j={b}_j-{b}_j^2\epsilon_j+\beta_j{b}_j^2\epsilon^2_j\nonumber\\
\beta_{j(k)}&=&{b}_{j(k)}-\beta_{j(k)}{b}_{j(k)}\epsilon_{j(k)}={b}_{j(k)}-{b}_{j(k)}^2\epsilon_{j(k)}+\beta_{j(k)}{b}_{j(k)}^2\epsilon^2_{j(k)}
\label{eq45}
\eqa
with ${b}_j=\frac{1}{1+E\bgma_j^*{\bf A}_j\bgma_j}$, $\epsilon_j=\bgma_j^*{\bf A}_j\bgma_j-E\bgma_j^*{\bf A}_j\bgma_j$, and $b_{j(k)}$ and $\epsilon_{j(k)}$ are similarly defined by replacing $\bbA_j^{-1}$ as $\bbA_{jk}^{-1}$.
By the same manner, we can decompose $\bgma_{j_2}^*\bbA\bgma_{k_2}$ into similar 6 terms and then we will estimate the expectations of the 36 products in the expansion
of $\bgma_{j_1}^*\bbA\bgma_{k_1}(\bgma_{j_2}^*\bbA\bgma_{k_2})^*$.

\noindent
{\bf Case 1. There are at least 6 terms of $\bbA_{j_1j_2k_1,k_2}^{-1}:=\bbB$'s contained in the product.} We shall use the fact that all $\beta$-factors are bounded $|z|/v\le K$. Then we can show that the term is bounded by $O(n^{-3})$. Say, for the product of the two 6-th terms, its expectation is bounded by
\bqn
&&KE\Big|(\bgma_{j_1}^*\bbB\bgma_{j_2}\bgma_{j_2}^*\bbB\bgma_{k_2}
\bgma_{k_2}^*\bbB\bgma_{j_2}\bgma_{j_2}^*\bbB\bgma_{k_1})(\bgma_{j_2}^*\bbB\bgma_{j_1}\bgma_{j_1}^*\bbB\bgma_{k_1}
\bgma_{k_1}^*\bbB\bgma_{j_1}\bgma_{j_1}^*\bbB\bgma_{k_2})^*\Big|\\
&\le&K\left(\rE\Big|(\bgma_{j_1}^*\bbB\bgma_{j_2}\bgma_{j_1}^*\bbB\bgma_{k_1}
\bgma_{k_2}^*\bbB\bgma_{j_2}\bgma_{j_2}^*\bbB\bgma_{k_1})\Big|^2\rE\Big|(\bgma_{j_2}^*\bbB\bgma_{j_1}\bgma_{j_2}^*\bbB\bgma_{k_2}
\bgma_{k_1}^*\bbB\bgma_{j_1}\bgma_{j_2}^*\bbB\bgma_{k_2})\Big|^2\right)^{1/2}.
\eqn
Note that the factors $\bgma_{j_2}^*\bbB\bgma_{k_2}$ in the first batch and $(\bgma_{j_1}^*\bbB\bgma_{k_1})^*$ in the second batch are exchanged positions
when using the Cauchy-Schwarz for avoiding 8the power of the $\bgma$ under the expectation sign.

Applying the formula
\bqn
&&\rE\left|\sum_{i=1}^n\bar X_iY_i\sum_{j=1}^n \bar X_jZ_j\sum_{k=1}^n \bar Y_kZ_k\right|^2\\
&=&\rE\left(\sum_{i\ne j}\left(|Y_i|^2|Z_j|^2+Y_i\bar Z_i\bar Y_j Z_j+|\rE X_i^2|^2Y_iZ_i\bar Y_j \bar Z_j\right)+\sum_{i=1}^n
\rE|X_i|^4|Y_i|^2|Z_i|^2\right)\left|\sum_{k=1}^n\bar Y_k Z_k\right|^2\\
&=&\left(\sum_{i\ne k}\left(2(n-2)|Y_i|^2|Y_k|^2+(n-2)(|\rE Z_i^2|^2+|\rE X_i^2|^2)Y_i^2\bar Y_j^2\right)+\sum_{i=1}^n
\rE|X_i|^4|Y_i|^4\rE|Z_i|^4\right)\\
&\le &2\kappa(n-2)\left(\sum_{i=1}|Y_i|^2\right)^2+\max_i\{\rE|X_i|^4\rE|Z_i|^4-2\kappa\}\sum_{i=1}^n|Y_i|^4,
\eqn
where $\kappa=2$ for the real case and 1 for the complex case, $X_i, Z_k$ are independent random variables with mean 0, variance 1 and finite 4th moment, and further
$\rE X_i^2=0$ (and $\rE Z_i^2=0$) if they are complex,
we will have
\bqn
&&\rE\Big|(\bgma_{j_1}^*\bbB\bgma_{j_2}\bgma_{j_1}^*\bbB\bgma_{k_1}
\bgma_{k_2}^*\bbB\bgma_{j_2}\bgma_{j_2}^*\bbB\bgma_{k_1})\Big|^2=\frac1n\rE|(\bgma_{j_1}^*\bbB\bgma_{j_2}\bgma_{j_1}^*\bbB\bgma_{k_1}\bgma_{j_2}^*\bbB\bgma_{k_1})\Big|^2
\bgma_{j_2}^*\bbB\bbT\bbB^*\bgma_{j_2}\\
&\le& \frac{K}{n^4}\rE(\bgma_{j_2}^*\bbB\bbT\bbB^*\bgma_{j_2})^3+\frac{K}{n^5}\sum_{i=1}^n\rE\Big|\bbe_i'\bbT^{1/2}\bbB\bgma_{j_2}\Big|^4\bgma_{j_2}^*\bbB\bbT\bbB^*\bgma_{j_2}\\
&\le&\frac{K}{n^4}\bigg[\left(\frac1n(\bbB\bbT\bbB^*\bbT)\right)^3+\frac1{n^3}\sum_{i=1}^n|\bbe_i'\bbT^{1/2}\bbB\bbT\bbB^*\bbT^{1/2}\bbe_i|^6\rE|X_{ij_2}|^6\bigg]
=O(n^{-4}),
\eqn
where $\bbe_i$ is the standard $i$-th unit $p$-vector, i.e., its $i$-th entry is 1 and other $p-1$ entries 0. In the last step of the above derivation, we have used facts that
$\rE|X_{ij_2}^6|\le \eta_n^2n\max \rE\rE|x_{ij}^4|=o(n)$ and $\bbe_i'\bbT^{1/2}\bbB\bbT\bbB^*\bbT^{1/2}\bbe_i\le \|\bbT\|^2/v^2$.

By similar approach, one can prove that the expectation of other products with the number of $\bbB$ less than or equal to 6 are bounded by $O(n^{-3})$.

\noindent
{\bf Case 2. There are 5 $\bbB$'s contained in the product.} We shall use the first expansion of $\beta_{j_1}$ and $\beta_{j_2}$ and then use the bound
 bounded $|z|/v\le K$ for $\beta$'s. Then we can show that such terms are also bounded by $O(n^{-3})$. Say, for the product of the first term of the first factor and the 6-th term
 of the second factor, its expectation is bounded by
\bqn
&&\bigg|\rE(\beta_{j_1}\beta_{k_1(j_1)}\bgma_{j_1}^*\bbB\bgma_{k_1})(\beta_{j_2}\beta_{k_2(j_2)}\beta_{k_1(j_2k_2)}\beta_{j_1(j_2k_1k_2)}^2\bgma_{j_2}^*\bbB\bgma_{j_1}\bgma_{j_1}^*\bbB\bgma_{k_1}
\bgma_{k_1}^*\bbB\bgma_{j_1}\bgma_{j_1}^*\bbB\bgma_{k_2})^*\Big|\\
&=&\bigg|\rE\Big(\beta_{j_1}\beta_{k_1(j_1)}\beta_{j_2}\beta_{k_2(j_2)}\beta_{k_1(j_2k_2)}\beta_{j_1(j_2k_1k_2)}^2
-b_{j_1}b_{k_1(j_1)}b_{j_2}1_{k_2(j_2)}b_{k_1(j_2k_2)}b_{j_1(j_2k_1k_2)}^2\Big)\times\\
&&\ \ \bgma_{j_1}^*\bbB\bgma_{k_1}\bgma_{j_2}^*\bbB\bgma_{j_1}\bgma_{j_1}^*\bbB\bgma_{k_1}
\bgma_{k_1}^*\bbB\bgma_{j_1}\bgma_{j_1}^*\bbB\bgma_{k_2})^*\Big|\\
&\le&K\bigg(\rE\Big|\Big(\beta_{j_1}\beta_{k_1(j_1)}\beta_{j_2}\beta_{k_2(j_2)}\beta_{k_1(j_2k_2)}\beta_{j_1(j_2k_1k_2)}^2
-b_{j_1}b_{k_1(j_1)}b_{j_2}1_{k_2(j_2)}b_{k_1(j_2k_2)}b_{j_1(j_2k_1k_2)}^2\Big)\\
&&\ \ (\bgma_{j_1}^*\bbB\bgma_{k_1}\bgma_{j_2}^*\bbB\bgma_{j_1})
\Big|^2\rE\Big|(\bgma_{j_1}^*\bbB\bgma_{k_1}\Big|^4\Big|
\bgma_{j_1}^*\bbB\bgma_{k_2}\Big|^2\bigg)^{1/2}\le O(n^{-3}).
\eqn
Here, we have used the fact that each term in the expansion of
$$\Big(\beta_{j_1}\beta_{k_1(j_1)}\beta_{j_2}\beta_{k_2(j_2)}\beta_{k_1(j_2k_2)}\beta_{j_1(j_2k_1k_2)}^2
-b_{j_1}b_{k_1(j_1)}b_{j_2}1_{k_2(j_2)}b_{k_1(j_2k_2)}b_{j_1(j_2k_1k_2)}^2\Big)$$
contains at leat one $\epsilon$ function which the centralized quadratic form of $\bgma$.
The use the same approach employed in Case1, one can show that the bound id $O(n^{-3})$.

\noindent
{\bf Case 3. There are less than 5 $\bbB$'s contained in the product.} If the number is 4, we need to further expand the matrix $\bbA_{j_1}$ in $\epsilon_{j_1}$
as $\bbA_{j_1j_2}^{-1}-\bbA_{j_1j_2}^{-1}\bgma_{j_2}\bgma_{j_2}^*\bbA_{j_1j_2}^{-1}\beta_{j_2(j_1)}$,
 expand $\bbA_{j_2}^{-1}=\bbA_{j_1j_2}^{-1}-\bbA_{j_1j_2}^{-1}\bgma_{j_1}\bgma_{j_1}^*\bbA_{j_1j_2}^{-1}\beta_{j_1(j_2)}$ in $\epsilon_{j_2}$,
 and then use the approach employed in Case 2 to obtain the desired bound.

 If the number is less than 4, we need to further expand the inverses of $\bbA$-matrices. The details are omitted. Finally, we obtain that
$$
\sum_{(4)}=O(\frac{1}{n}).
$$

Similarly, we have $$
\sum_{(3)}=O(\frac{1}{n}).
$$
Because $\rtr(\bbA^{-2}\gD)=\frac{d}{dz}\rtr(\bbA^{-1}\gD)$,
then we have
$$
\rE|\rtr(\bbA^{-2}\gD)|^2=\rE\left|\frac1n\sum_{j\ne k}\bgma_j^*\bbA^{-2}\bgma_k\right|^2=O(\frac{1}{n}).
$$
The lemma is proved.

\begin{lem}\label{l3}After truncation and normalization, we have
$$
\rtr(\bbA^{-2}\gD\bbA^{-1}\gD)=(\underline{m}_{y}(z)+z\underline{m}_{y}'(z))(1+z\underline{m}_{y}(z))
$$
in $L_2$ uniformally for $z\in{\mathbb C}^+$.
\end{lem}

\proof \ Set
$\rtr\bbA^{-1}(z_1)\gD\bbA^{-1}(z_2)\gD=\frac1{n^2}\sum_{j\ne k, i\ne t}
\bgma_i^*\bbA^{-1}(z_1)\bgma_k\bgma_j^*\bbA^{-1}(z_2)\bgma_t
=Q_1+Q_2$ where
$$Q_1=\frac1{n^2}\sum_{j\not=k}^{n}\bgma_j^*\bbA^{-1}(z_1)\bgma_j\bgma_k^*\bbA^{-1}(z_2)\bgma_k
\quad\mbox{and}\quad
Q_2=\frac1{n^2}\sum_{j\ne k, i\ne t\atop
i\ne k,{\rm or} j\ne t}\bgma_i^*\bbA^{-1}(z_1)\bgma_k\bgma_j^*\bbA^{-1}(z_2)\bgma_t.
$$
By Lemma \ref{l1} and \ref{l2}, we obtain
$\rE|Q_1-(1+z\underline{m}_{y}(z_1))(1+z\underline{m}_{y}(z_2))|^2\leq Kn^{-1}$ and
$\rE|Q_2|^2=o(1).$
We thus have
$\rE|\rtr\bbA^{-1}(z_1)\gD\bbA^{-1}(z_2)\gD-(1+z\underline{m}_{y}(z_1))(1+z\underline{m}_{y}(z_2))|^2=o(1).$
Consequently, because $\rtr\bbA^{-2}(z_1)\gD\bbA^{-1}(z_2)\gD=\frac{\partial rtr\bbA^{-1}(z_1)\gD\bbA^{-1}(z_2)\gD}{\partial z_1}$, then we have
$
\rE|\rtr\bbA^{-2}(z_1)\gD\bbA^{-1}(z_2)\gD-\frac{\partial}{\partial z_1}g(z_1)g(z_2)|^2=o(1).
$
That is,
$$
\rtr\bbA^{-2}(z_1)\gD\bbA^{-1}(z_2)\gD=g(z_2)g'(z_1) \mbox{in $L_2$.}
$$
By setting $z_1=z_2=z$, we obtain
$\rtr(\bbA^{-2}\gD\bbA^{-1}\gD)=g(z)g'(z) \mbox{in $L_2$.}$
\begin{lem}\label{l4} After truncation and normalization, we have
$\rtr(\bbA^{-1}\gD)^3(\bbA-\gD)^{-1}=g(z)\rtr((\bbA^{-1}\gD)^2(\bbA-\gD)^{-1})
+o_p(1)$
uniformly for $z\in{\mathbb C}^+$.
\end{lem}
\proof\ We have
\bqn
&&\rtr(\bbA^{-1}\gD)^3(\bbA-\gD)^{-1}
=E\frac1{n^3}\sum_{i\ne t, j\ne g\atop h\ne s}\bgma_i^*\bbA^{-1}\bgma_j\bgma_g^*\bbA^{-1}\bgma_h\bgma_s^*
(\bbA-\gD)^{-1}\bbA^{-1}\bgma_t\nonumber\\
&=&\frac1{n^3}\sum_{i\ne t, j\ne g\atop i=j, h\ne s}\bgma_i^*\bbA^{-1}\bgma_j\bgma_g^*\bbA^{-1}\bgma_h
\bgma_s^*(\bbA-\gD)^{-1}\bbA^{-1}\bgma_t
+\frac1{n^3}\sum_{i\ne t, j\ne g\atop i\ne j, h\ne s}\bgma_i^*\bbA^{-1}\bgma_j\bgma_g^*\bbA^{-1}\bgma_h
\bgma_s^*(\bbA-\gD)^{-1}\bbA^{-1}\bgma_t\nonumber\\
&=&g(z)\frac1{n^2}\sum_{h\ne s}\bgma_g^*\bbA^{-1}\bgma_h\bgma_s^*(\bbA-\gD)^{-1}\bbA^{-1}\bgma_t
+o_p(1)\\
&=&g(z)\rtr((\bbA^{-1}\gD)^2(\bbA-\gD)^{-1})+g(z)\frac1{n^2}\sum_{g=t\atop h\ne s}\bgma_g^*\bbA^{-1}
\bgma_h\bgma_s^*(\bbA-\gD)^{-1}\bbA^{-1}\bgma_t
+o_p(1)\\
&=&g(z)\rtr((\bbA^{-1}\gD)^2(\bbA-\gD)^{-1})
+o_p(1).
\eqn
Then by Lemma \ref{l3}, we have
\bqn
\rtr(\bbA^{-1}\gD)^2(\bbA-\gD)^{-1}&=&\rtr(\bbA^{-1}\gD)^2(\bbA)^{-1}+\rtr(\bbA^{-1}\gD)^3(\bbA-\gD)^{-1}\\
&=&\rtr(\bbA^{-1}\gD)^2(\bbA)^{-1}+g(z)\rtr(\bbA^{-1}\gD)^2(\bbA-\gD)^{-1}+o_p(1)\\
&=&\frac{(1+z\underline{m}_{y}(z))(\underline{m}_{y}(z)+z\underline{m}_{y}'(z))}{1-g(z)}+o_p(1).
\eqn
Hence, we obtain the following lemma.
\begin{lem}\label{thm2}After truncation and normalization, we have
$$
tr{\bf A}^{-2}(z)\gD
+tr{\bf A}^{-1}(z)(\gD{\bf A}^{-1}(z))^2+tr\left({\bf A}(z)-\gD\right)^{-1}(\gD{\bf A}^{-1}(z))^3=
\frac{(\underline{m}_{y}(z)+z\underline{m}_{y}'(z))(1+z\underline{m}_{y}(z))}{-z\underline{m}_{y}(z)}+o_p(1)
$$
uniformly for $z\in{\mathbb C}^+$.
\end{lem}

\end{document}